\documentclass[reqno, a4paper]{amsart}

\usepackage{amssymb,amsmath,amsthm}
\newtheorem{thm}{Theorem}
\newtheorem{lem}[thm]{Lemma}
\newtheorem{cor}[thm]{Corollary}
\newtheorem{rem}[thm]{Remark}
\newtheorem{prop}[thm]{Proposition}
\newtheorem{dfn}[thm]{Definition}

\begin{document}

\title[Traces generated by exponentiation invariant generalised limits]{Dixmier traces generated by exponentiation invariant generalised limits}
\author{Fedor Sukochev}
\email{f.sukochev@unsw.edu.au}
\author{Alexandr Usachev}
\email{a.usachev@unsw.edu.au (Corresponding Author)}
\author{Dmitriy Zanin}
\email{d.zanin@unsw.edu.au}
 \address{School of Mathematics and Statistics, University of New South Wales, Sydney, 2052, Australia.}
\maketitle

\begin{abstract}
We define a new class of singular positive traces on the ideal $\mathcal M_{1,\infty}$ of $B(H)$ generated by exponentiation invariant generalized limits. 
We prove that this new class is strictly contained in the class of all Dixmier traces. We also prove a Lidskii-type formula for this class of traces.
\end{abstract}



 \section{Introduction and Preliminaries}

In the framework of noncommutative geometry the Dixmier traces, originally introduced by J.~Dixmier in~\cite{D}, have become  an 
indispensable tool. These traces are defined via dilation invariant generalized limits $\omega$ on $L_\infty(0,\infty)$.
Depending  on a specific problem, various additional conditions on $\omega$ are imposed~\cite{BF}-\cite{D},~\cite{LSS},~\cite{Sed}-\cite{SUZ2}.

For instance, important formulae of noncommutative geometry, involving heat kernel asymptotics and generalised $\zeta-$functions 
(see e.g.~\cite{CPS2,CPRS1}) were established for Dixmier traces ${\rm Tr}_\omega$, 
provided that $\omega$ was additionally chosen to be exponentiation invariant.
Indeed, in~\cite{CPS2} these formulae were proven for Dixmier traces, 
generated by Ces\`aro and exponentiation invariant generalized limits $\omega$.
In~\cite{CPRS1} these assumptions were relaxed to dilation and exponentiation invariance (see the definition below).

We thank the referee for very careful reading of the manuscript and suggesting a number of improvements. 
In particular, the statement of Theorem~\ref{zeta} is due to the referee.

\subsection{Generalized limits}

Let $L_\infty=L_\infty(0,\infty)$ be the space of all real-valued
bounded Lebesgue measurable functions on $(0,\infty)$ equipped with the norm
$$\|x\|_{L_\infty}:=\mathop{\rm esssup}\limits_{t>0} |x(t)|.$$

A normalized positive linear functional on $L_\infty$ which equals the ordinary limit on convergent (at infinity) sequences is called a generalized limit.

Define the following linear transformations on $L_\infty$

1. Translations  $$(T_lx)(t) := x(t+l), \ l>0, \ x\in L_\infty;$$

2. Dilations  $$(\sigma_{\frac1\beta} x)(t):=x(\beta t), \ \beta>0, \ x\in L_\infty;$$

3. Exponentiations $$(P^ax)(t)=x(t^a), \ a>0, \ x\in L_\infty.$$

A generalized limit $\omega$ on $L_\infty$ is said to be dilation invariant if 
$$\omega(\sigma_s x)=\omega(x) \ \text{for every} \ x\in L_\infty \ \text{and every} \ s>0.$$ 

Similarly we define translation and exponentiation invariant generalized limits.

Translation and dilation invariant generalized limits on $L_\infty$ were studied in~\cite{LP,Sed, SUZ1} in terms of Banach-type sublinear functionals.

It was proven in~\cite[Theorem 13]{SUZ1} (see also~\cite{LP,Sed}) that 
for every translation invariant generalized limit $\omega$ on $L_\infty$ and any uniformly continuous function $x\in L_\infty$ 
we have
\begin{equation}\label{sucheston}
\omega(x) \le p_T(x):= \lim_{t\to+\infty}\sup_{h\ge0}\frac1t\int_0^t x(s+h) ds. 
\end{equation}
It should be pointed out that the condition that $x$ is uniformly continuous is essential (see e.g.~\cite[Remark 5.6 (1)]{LP}).

\begin{rem}\label{exp}
 For every translation invariant generalized limit $\gamma$ on $L_\infty$ 
the composition $\gamma \circ \exp$ defines dilation invariant generalized limit (see~\cite[Remark 16]{SUZ1} for details).
Similarly, $\gamma \circ \exp \circ \exp$ defines an exponentiation invariant generalized limit.

Conversely, if $\omega$ is an exponentiation invariant generalized limit,
then $\omega\circ\log$ and $\omega\circ\log\circ\log$ are dilation and translation invariant generalized limits, respectively.
\end{rem}

Using Remark~\ref{exp} the result similar to~\eqref{sucheston} was proven for dilation invariant generalized limit and 
any function $x\in L_\infty$ such that $x \circ \exp$ is uniformly continuous~\cite[Theorem 17]{SUZ1}.

\subsection{Dixmier traces}
Let $B(H)$ be an algebra of all bounded linear operators on a separable Hilbert space
$H$ equipped with the uniform norm and let ${\rm Tr}$ be the standard trace.

For every operator $T\in B(H)$ a generalized singular value function $\mu(T)$ is defined by the formula
$$\mu(t, T)=\inf \{\|Tp\| : \ p \ \text{is a projection in} \ B(H) \ \text{with} \ {\rm Tr}(1-p)\le t\}.$$
For a compact operator $T$, it can be proven that $\mu(k-1, T)$ is the $k$-th largest eigenvalue of an operator $|T|=(T^*T)^{1/2}$, $k \in \mathbb N$.

Recall that for any $T\in B(H)$ the distribution function of $T$ is defined (see~\cite{FK}) by setting
$$
d_T(t):= {\rm Tr} (\chi_{(t,\infty)}(|T|)),\quad t>0.
$$
Here $\chi_{(t,\infty)}(|T|)$ is the spectral projection of $|T|$
corresponding to the interval $(t,\infty)$. 
By~\cite[Proposition 2.2]{FK}, we have
$$\mu(s,T)=\inf\{t\ge 0\ :\ d_T(t)\leq s\},
$$
so, we infer that for any operator $T$, the distribution
function $d_T$ coincides with the (classical) distribution
function of $\mu(\cdot,T)$.

Since $B(H)$ is an atomic von Neumann algebra and traces of all atoms equal to 1, it follows that 
$\mu(T)$ is a step function and 
$$\mu(T)=\sum_{n=0}^\infty \mu(n,T) \chi_{[n,n+1)}\ \text{for every} \ T\in B(H).$$ 

The classical Dixmier ideal (see e.g.~\cite{CS,KSS,LSS}) is defined as
$$\mathcal M_{1,\infty}:= \left\{T: \|T\|_{\mathcal M_{1,\infty}}:= \sup_{t> 0}  \frac{1}{\log(1+t)} \int_0^t \mu(s,T) \ ds < \infty\right\}.$$
A definition of (Dixmier) traces given in~\cite{D} can be restated as follows
\begin{equation}\label{Dix_dfn}
{\rm Tr}_\omega (T):= \omega \left( \frac{1}{\log(1+t)} \int_0^t \mu(s,T) \ ds\right), \quad 0\le T\in \mathcal M_{1,\infty}, 
\end{equation}
where $\omega$ is an arbitrary dilation invariant generalized limit on $L_\infty$. We denote the set of all Dixmier traces by $\mathcal D$.

\subsection{Measurability}
The following natural way to generate dilation invariant generalized limits was suggested in~\cite[Section IV, 2$\beta$]{C}.  A.~Connes 
observed that for any generalised limit $\gamma$ on $L_\infty$ 
a functional $\omega:=\gamma\circ M$ is a dilation invariant generalized limit on $L_\infty$.
Here, the operator $M:L_\infty\to L_\infty$ is the Ces\`aro mean of the multiplicative group $\mathbb R_+$, given by the formula
$$(Mx)(t):=\frac1{\log t} \int_1^t x(s)\,\frac{ds}s.$$
The subclass $\mathcal{C}\subset \mathcal{D}$ of all Dixmier traces ${\rm Tr}_{\omega}$ defined by such $\omega$'s 
was termed Connes-Dixmier traces in~\cite{LSS}. It was proven in~\cite{SUZ2} that $\mathcal{C}\subsetneqq \mathcal{D}$.

As it was mentioned above, various important formulae of noncommutative geometry were established for 
dilation and exponentiation invariant generalized limits $\omega$.
The former assumption was needed in order that the formula~\eqref{Dix_dfn} defines a Dixmier trace.

In the present paper we prove that for every exponentiation invariant generalized limit $\omega$ on $L_\infty$ the formula~\eqref{Dix_dfn} defines a Dixmier trace (see Proposition~\ref{linearity lemma} below). Denote by $\mathcal D_P$ the set of all Dixmier traces generated by such $\omega$'s.

The following definitions were motivated by A.~Connes~\cite[IV.2.$\beta$.Definition 7]{C} (see also~\cite[Definition 3.2]{LSS}).
\begin{dfn}\label{Dmeas}
An operator $T\in \mathcal M_{1,\infty}$ is called Dixmier measurable if ${\rm Tr}_\omega(T)$ takes the same value for all ${\rm Tr}_\omega\in \mathcal{D}$.
\end{dfn}

A criterion for an operator 
$T\in \mathcal M_{1,\infty}$ (respectively, a positive operator $T\in \mathcal M_{1,\infty}$) 
to be Dixmier measurable can be found in~\cite[Theorem 21]{SUZ1} (respectively,~\cite{LSS}).
The crucial point in these proof is played by the fact that the function 
$$t \to \frac{1}{\log(1+e^t)} \int_0^{e^t} \mu(s,T) \ ds$$ is uniformly continuous.

\begin{dfn}\label{Pmeas}
An operator $T\in \mathcal M_{1,\infty}$ is called $\mathcal D_P$-measurable if ${\rm Tr}_\omega(T)$ takes the same value for all ${\rm Tr}_\omega\in \mathcal D_P$.
\end{dfn} 

Using Remark~\ref{exp} we can easily write down the Banach-type sublinear functional for exponentiation invariant generalized limits.
However, we cannot gain any results about $\mathcal D_P$-measurability on this way, 
since there exist $0\le T\in \mathcal M_{1,\infty}$ such that the function 
$$t \to \frac{1}{\log(1+e^{e^t})} \int_0^{e^{e^t}} \mu(s,T) \ ds$$ 
is not uniformly continuous. 
As an example of such $T$ we may take the operator $T_0$ defined in Theorem~\ref{PD_ex} below.

We prove (see Theorem~\ref{PD_ex} below)  that the class of 
$\mathcal D_P$-measurable operators is stricly wider then the class of Dixmier measurable operators.
In particular, $\mathcal D_P$ is a proper subset of $\mathcal D$.

\subsection{Lidskii formula for Dixmier traces}
The classical Lidskii Theorem asserts that
$${\rm Tr}(T)=\sum_{n\geq 0}\lambda(n,T)$$
for any trace class operator $T.$ Here, $\{\lambda(n,T)\}_{n\geq 0}$ is the sequence of eigenvalues of $T$ 
(counting with algebraic multiplicities), taken in an arbitrary order.
This arbitrariness of the order is due to the absolute convergence of the series $\sum_{n\geq 0}|\lambda(n,T)|.$ 
In particular, we can choose a decreasing order of $|\lambda(n,T)|$.

The core difference of this situation with the setting of Dixmier traces living on the ideal $\mathcal{M}_{1,\infty}$ 
is that the series $\sum_{n\geq 0} |\lambda(n,T)|$ 
diverges for $T\in\mathcal{M}_{1,\infty}$ (see~\cite{SSZ} for a detailed explanation).

The following theorem gives an analogue of Lidskii formula for Connes-Dixmier traces.
\begin{thm}~\cite[Theorem 2]{SSZ}\label{formulalidskogo} Let ${\rm Tr}_{\omega}$ be a Connes-Dixmier trace on $\mathcal{M}_{1,\infty}.$ We have
\begin{equation}\label{as lidskiy formula}
{\rm Tr}_\omega(T)=\omega\left(\frac1{\log(1+t)}\sum_{\lambda\in\sigma(T): |\lambda|>1/t}\lambda\right),\quad T\in\mathcal{M}_{1,\infty},
\end{equation}
where $\sigma(T)$ denotes the spectrum of an operator $T$.
\end{thm}

It was also shown in~\cite[Theorem 5]{SSZ} that there exists a Dixmier trace ${\rm Tr}_\omega$ 
such that the formula~\eqref{as lidskiy formula} does not hold. 

One of the main results of the present paper (see Theorem~\ref{LFP}) asserts that the formula~\eqref{as lidskiy formula} holds for every Dixmier trace
${\rm Tr}_\omega \in {\mathcal D}_P$.

 \section{Dixmier traces generated by exponentiation invariant generalised limits}
The following proposition is an analogue of~\cite[Proposition 10]{KSS} for exponentiation invariant generalised limits. 
It shows that every exponentiation invariant generalised limit generates a Dixmier trace.

\begin{prop}\label{linearity lemma} 
For every exponentiation invariant generalised limit $\omega$ on $L_\infty$ the weight 
$$
{\rm Tr}_\omega(T):=\omega \left( \frac{1}{\log(1+t)} \int_0^t\mu(s,T) \ ds\right), \quad 0\le T\in \mathcal M_{1,\infty},
$$ 
extends to a non-normal trace on $\mathcal M_{1,\infty}$.
\end{prop}

\begin{proof} By the construction of ${\rm Tr}_\omega$ we only need to prove its additivity on the positive cone of $\mathcal M_{1,\infty}$.

Let $0 \le A,B\in \mathcal M_{1,\infty}$. By~\cite[Theorem 4.4 (ii)]{FK} we have 
$$\int_0^{t/2}\mu(s,A)+\mu(s,B)\,ds \le \int_0^t \mu(s,A+B)\, ds \le \int_0^t \mu(s,A)+\mu(s,B)\, ds $$ 
for every $t>0$.

Thus, using the positivity of $\omega$, we obtain
$${\rm Tr}_{\omega}(A+B) \le {\rm Tr}_{\omega}(A)+{\rm Tr}_{\omega}(B)$$
for every generalised limit $\omega$.

On the other hand, we have
\begin{align*}
{\rm Tr}_{\omega}(A)+{\rm Tr}_{\omega}(B)&=\omega\left( \frac1{\log(1+t)}\int_0^{t} \mu(s,A)+\mu(s,B)\,ds\right)\\ 
&\leq\omega\left( \frac1{\log(1+t)}\int_0^{2t}\mu(s,A+B) ds\right).
\end{align*}

Since for every $\varepsilon>0$ we have $2t\le t^{1+\varepsilon}$ (for $t$ large enough) and 
$$\log(1+t^{1+\varepsilon})\le(1+\varepsilon)\log(1+t),$$
the following estimate holds
\begin{align*}
{\rm Tr}_{\omega}(A)+{\rm Tr}_{\omega}(B)&\leq\omega\left( \frac1{\log(1+t)}\int_0^{2t}\mu(s,A+B)ds\right)\\
&\leq (1+\varepsilon) \, \omega\left(\frac1{\log(1+t^{1+\varepsilon})} \int_0^{t^{1+\varepsilon}} \mu(s,A+B)ds\right)\\
&\leq (1+\varepsilon) \, \omega\left(\frac1{\log(1+t)} \int_0^{t} \mu(s,A+B)ds\right)\\
&=(1+\varepsilon) \, {\rm Tr}_{\omega}(A+B).
\end{align*}

Since $\varepsilon>0$ can be chosen arbitrary small, we conclude that ${\rm Tr}_{\omega}(A)+{\rm Tr}_{\omega}(B) \leq {\rm Tr}_{\omega}(A+B)$.

\end{proof}
It is obvious from the definition of ${\rm Tr}_\omega$ in Proposition~\ref{linearity lemma} that every such functional is fully symmetric 
(see e.g.~\cite{CS, KSS}). By the main result of~\cite{KSS} the set of all fully symmetric functionals on $\mathcal M_{1,\infty}$ coincides with
the set of all Dixmier traces, that is with the set ${\mathcal D}$.
Thus, singular traces generated by an exponentiation invariant generalised limit $\omega$ are Dixmier traces.
Next, we show that the class of Dixmier measurable operators is strictly wider then the class of ${\mathcal D}_P$-measurable operators.

First, we state two auxiliary lemmas.

\begin{lem}\label{periodicity}\cite[Example 5.6 (1)]{LP}
 Let $x\in L_\infty$ be a locally Riemann integrable function. 
If $x$ is a periodic function and its period is $l>0$, then
$$\gamma(x) = \frac1l \int_0^l x(s) ds,$$
for every translation invariant generalized limit $\gamma$ on $L_\infty$.
\end{lem}

The following lemma is an analogue of~\cite[Lemma 1.2]{SUZ2} for exponentiation invariant generalised limits.
\begin{lem}\label{aux0}
For every $T\in \mathcal M_{1,\infty}$ and for every exponentiation invariant generalised limit $\omega$, we have
$$\omega\left( \frac{t \mu(t,T)}{\log(1+t)} \right) =0.$$
\end{lem}

\begin{proof}
 Since $\omega$ is an exponentiation invariant generalised limit, it follows that
$$\omega\left( \frac1{\log(1+t)} \int_0^{t/2} \mu(s,T) \ ds \right) = \omega\left( \frac1{\log(1+t^{1+\varepsilon})} \int_0^\frac{t^{1+\varepsilon}}{2} \mu(s,T) \ ds \right)$$
for every $\varepsilon>0$.
For every fixed $\varepsilon>0$ we have $\frac{t^{1+\varepsilon}}{2} \ge t$ (for $t$ large enough) and
$$\log(1+t^{1+\varepsilon})\le(1+\varepsilon)\log(1+t), \ \forall t>0.$$

Hence, 
$$\omega\left( \frac1{\log(1+t)} \int_0^{t/2} \mu(s,T) \ ds \right) 
\ge \frac1{1+\varepsilon} \cdot \omega\left( \frac1{\log(1+t)} \int_0^t \mu(s,T) \ ds \right).$$

Since $\varepsilon$ is abritrary, it follows that
$$\omega\left( \frac1{\log(1+t)} \int_0^{t/2} \mu(s,T) \ ds \right) =\omega\left( \frac1{\log(1+t)} \int_0^t \mu(s,T) \ ds \right).$$

Consequently, since $\mu(T)$ is decreasing, it follows that
$$0=\omega\left( \frac1{\log(1+t)} \int_{t/2}^t \mu(s,T) \ ds \right) \ge \omega\left( \frac{t \mu(t,T)}{2\log(1+t)} \right)\ge 0.$$
\end{proof}

In view of the main result of~\cite{LSS} the set of all positive ${\mathcal D}$-measurable and the set of all positive ${\mathcal C}$-measurable
elements coincide. The following theorem shows that the sets of all positive ${\mathcal D}_P$-measurable  
and all positive ${\mathcal C}$-measurable elements are distinct, in particular, the classes ${\mathcal D}_P$ and ${\mathcal C}$ are different.

\begin{thm}\label{PD_ex}
There exists a positive Dixmier non-measurable operator $T_0\in \mathcal M_{1,\infty}$, 
such that all ${\rm Tr}_\omega \in {\mathcal D}_P$ take the same value on $T_0$.
\end{thm}
\begin{proof}
Let $T_0$ be a compact operator on $H$ such that
$$\mu(T_0)=\sup_{k\geq0}e^{k-e^k}\chi_{[0,\lfloor e^{e^k}\rfloor)}.$$ 
Consider a function
$$z=\sup_{k\geq0}e^{k-e^k}\chi_{[0,e^{e^k})}.$$ 
It is easy to see that $\mu(T_0)-z \in L_1(0,\infty)\cap L_\infty(0,\infty)$.

For every $e^{e^n}\le t< e^{e^{n+1}}$ we have
 \begin{equation}\label{a1}
 \begin{aligned}
\frac{1}{\log(1+t)} \int_0^t\mu(s,T_0) \ ds&=\frac1{\log(1+t)}\left(\int_0^{e^{e^n}} z(s)ds+(t-e^{e^n})\mu(t, T_0)+O(1)\right)\\
&=\frac1{\log(1+t)}\left(\sum_{k=1}^{n}\int_{e^{e^{k-1}}}^{e^{e^k}}e^{k-e^k}ds+t\mu(t, T_0)+O(1)\right)\\
&= \frac1{e-1}\frac{e^{n+1}}{\log t}+\frac{t\mu(t, T_0)}{\log (1+t)}+o(1).
 \end{aligned} 
 \end{equation}

Now, it is easy to check that 
$$\frac{1}{\log(1+t)} \int_0^t\mu(s,T_0) \ ds \le \frac{e^2}{e-1} +o(1), \ t>0,$$ 
and, therefore $T_0\in \mathcal M_{1,\infty}$.
However, the limit
$$\lim_{t\to\infty} \frac{1}{\log(1+t)} \int_0^t\mu(s,T_0) \ ds$$
does not exist, since by~\eqref{a1}
$$\lim_{n\to\infty} \frac{1}{\log(1+e^{e^n})} \int_0^{e^{e^n}}\mu(s,T_0) \ ds = \frac{e}{e-1}$$
and
$$\lim_{n\to\infty} \frac{1}{\log(1+\frac{e^{e^{n+1}}+e^{e^n}}2)} \int_0^{\frac{e^{e^{n+1}}+e^{e^n}}2}\mu(s,T_0) \ ds = \frac1{e-1}+\frac12.$$

Thus,  $T_0$ is not Dixmier measurable (due to~\cite[Theorem 6.6]{LSS}).
It remains to show that $T_0$ is ${\mathcal D}_P$-measurable.

By Lemma~\ref{aux0} we have 
$$\omega\left( \frac{t \mu(t,T_0)}{\log(1+t)} \right) =0.$$

Denote 
$$x(t):=\frac{e}{e-1} \sum_{k=0}^\infty \frac{e^{k}}{\log t} \chi_{[e^{e^k}, e^{e^{k+1}})}(t) , \ t>0.$$

In view of~\eqref{a1} we only need to prove that all exponentiation invariant generalised limits coincide on $x$. 
By Remark~\ref{exp}, it is sufficient to show that all translation invariant generalised limits coincide on $x\circ \exp \circ \exp$.

Since for $n \le t< n+1$,
$$(x\circ \exp \circ \exp)(t) =\frac{e}{e-1} \sum_{k=0}^\infty\frac{e^k}{e^t}\chi_{[k,k+1)}(t),$$
we conclude that the function $x\circ \exp \circ \exp$ is a Riemann integrable periodic function.
By Lemma~\ref{periodicity}, for every translation invariant generalised limit $\gamma$ we have
$$\gamma(x\circ \exp \circ \exp) =\frac{e}{e-1} \int_0^1 \frac{ds}{e^s}=\frac{e}{e-1}\frac{e-1}{e}=1.$$
So, by Remark~\ref{exp}, $\omega(x)=1$ for every exponentiation invariant generalised limit $\omega$.

Hence,
$${\rm Tr}_\omega(T_0) = \omega\left(\frac{1}{\log(1+t)} \int_0^t\mu(s,T_0) \ ds\right)=1$$
for every exponentiation invariant generalised limit $\omega$.
\end{proof}

\begin{cor}
The set ${\mathcal D}_P$ is a proper subset of the set of all Dixmier traces.
\end{cor}

\section{Lidskii formula for Dixmier traces generated by exponentiation invariant generalised limits}

In the present section we first prove the Lidskii formula for positive operators $T\in \mathcal M_{1,\infty}$, 
then, using Ringrose's representation~\cite[Theorems 1,6,7]{R} of compact operators, we extend the formula to an arbitrary $T\in \mathcal M_{1,\infty}$.

\begin{lem}
 For every positive $T\in \mathcal M_{1,\infty}$ and for every exponentiation invariant generalized limit $\omega$ the following formula holds
$${\rm Tr}_\omega(T)=\omega\left(\frac1{\log(1+t)} \int_0^{d_T(1/t)} \mu(s,T) ds\right).$$
\end{lem}
\begin{proof}
Let $T$ be a positive operator from $\mathcal M_{1,\infty}$ and $\omega$ be an exponentiation invariant generalized limit on $L_\infty$.
 From the proof of \cite[Proposition 2.4]{CPS2} we know that
$$\int_0^t \mu(s,T) ds \le \int_0^{d_T(1/t)} \mu(s,T) ds +1.$$
Dividing both sides by $\log(1+t)$ and applying $\omega$, we obtain
$${\rm Tr}_\omega(T) \le \omega\left(\frac1{\log(1+t)} \int_0^{d_T(1/t)} \mu(s,T) ds\right).$$

On the other hand, by~\cite[Lemma 2.3]{CPS2} there exists a constant $C>0$ such that we eventually have
\begin{equation}\label{eq1}
d_T(\frac1t) \le C \cdot t \log (1+t)\le t^{1+\varepsilon},\ \varepsilon>0. 
\end{equation}

We have
\begin{equation}\label{eq2}
\frac1{\log(1+t)} = \frac{1+\varepsilon}{\log(1+t)^{1+\varepsilon}}\le \frac{1+\varepsilon}{\log(1+t^{1+\varepsilon})}. 
\end{equation}

Hence, by~\eqref{eq1} and~\eqref{eq2}, we obtain
\begin{align*}
\omega\left(\frac1{\log(1+t)} \int_0^{d_T(1/t)} \mu(s,T) ds\right) 
&\le \omega\left(\frac{1+\varepsilon}{\log(1+t^{1+\varepsilon})} \int_0^{t^{1+\varepsilon}} \mu(s,T) ds\right) \\
&\le (1+\varepsilon) \tau_\omega(T),
\end{align*}
since $\omega$ is an exponentiation invariant generalized limit.
Since $\varepsilon>0$ is arbitrary, we obtain the converse inequality.
\end{proof}

\begin{cor}\label{c1}
 If $T=T^*\in \mathcal M_{1,\infty}$, then
$${\rm Tr}_\omega(T)=\omega\left(\frac1{\log(1+t)} \sum_{\genfrac{}{}{0pt}{}{\lambda \in \sigma(T),}{|\lambda| >1/t}} \lambda\right)$$
for every exponentiation invariant Dixmier trace ${\rm Tr}_\omega$.
\end{cor}

The following Lemma is an analogue of~\cite[Lemma 42]{SUZ1} for exponentiation invariant generalized limits.

\begin{lem}\label{l1}
 For every positive $T\in \mathcal M_{1,\infty}$ and for every exponentiation invariant generalized limit $\omega$, we have
$$\omega\left(\frac1{t\log(1+t)} d_T(\frac1t)\right)=0.$$ 
\end{lem}
\begin{proof} 
Fix $0<\varepsilon<1$.
 By the definition of a distribution function we have
$$\frac1t \left( d_T(\frac1t)-d_T(\frac1{t^{1-\varepsilon}})\right) \le \sum_{1/t < \lambda \le 1/t^{1-\varepsilon}} \lambda.$$ 

Dividing both sides by $\log(1+t)$ and applying $\omega$, we obtain
$$\omega\left(\frac1{t\log(1+t)}  d_T(\frac1t)\right)-\omega\left(\frac1{t\log(1+t)}  d_T(\frac1{t^{1-\varepsilon}})\right) 
\le \omega\left(\frac1{\log(1+t)}\sum_{1/t < \lambda \le 1/t^{1-\varepsilon}} \lambda\right).$$ 
However, using~\eqref{eq1} we have
$$
\lim_{t\to\infty} \frac1{t\log(1+t)}  d_T(\frac1{t^{1-\varepsilon}}) \le 
\lim_{t\to\infty} \frac1{t\log(1+t)}  t^{(1-\varepsilon)(1+\varepsilon)}=0.
$$
Since $\omega$ is a generalized limit, it follows that 
$$\omega\left(\frac1{t\log(1+t)}  d_T(\frac1{t^{1-\varepsilon}})\right)=0. $$
Therefore, 
\begin{align*}
\omega\left(\frac1{t\log(1+t)}  d_T(\frac1t)\right) &\le \omega\left(\frac1{\log(1+t)}\sum_{1/t < \lambda \le 1/t^{1-\varepsilon}} \lambda\right)\\
&=\omega\left(\frac1{\log(1+t)}\sum_{\lambda > 1/t} \lambda\right)-
\omega\left(\frac1{\log(1+t)}\sum_{\lambda > 1/t^{1-\varepsilon}} \lambda\right)\\
&=\omega\left(\frac1{\log(1+t)}\sum_{\lambda > 1/t} \lambda\right)-
\omega\left(\frac{1-\varepsilon}{\log(1+t^{1-\varepsilon})}\sum_{\lambda > 1/t^{1-\varepsilon}} \lambda\right),
\end{align*}
where the last equality holds since $$\lim_{t\to\infty}\frac{(1-\varepsilon)\log(1+t)}{\log(1+t^{1-\varepsilon})}=1$$
and since $\omega$ is a generalized limit.

Using exponentiation invariance of $\omega$ and Corollary~\ref{c1}, we obtain
$$\omega\left(\frac1{t\log(1+t)}  d_T(\frac1t)\right) \le {\rm Tr}_\omega(T)-(1-\varepsilon){\rm Tr}_\omega(T)=\varepsilon \cdot {\rm Tr}_\omega(T).$$
Since $0<\varepsilon<1$ can be chosen arbitrarily small, we conclude that
$$\omega\left(\frac1{t\log(1+t)} d_T(\frac1t)\right)=0.$$
\end{proof}

We need two auxiliary lemmas.

\begin{lem}\label{l2}
 For every normal operator $T\in \mathcal M_{1,\infty}$ and for every exponentiation invariant generalized limit $\omega$, we have
$$\omega\left(\frac1{\log(1+t)} \sum_{\genfrac{}{}{0pt}{}{\lambda \in \sigma(T):}{ |\lambda| >1/t}} \lambda\right)=
\omega\left(\frac1{\log(1+t)} \sum_{\genfrac{}{}{0pt}{}{\lambda \in \sigma(T): |\Re( \lambda)| >1/t}{ \ or \ |\Im (\lambda)| >1/t}} \lambda\right).$$  
\end{lem}
\begin{proof} Consider the difference
\begin{align*}
\left| \sum_{\genfrac{}{}{0pt}{}{\lambda \in \sigma(T):}{ |\lambda| >1/t}} \lambda-
\sum_{\genfrac{}{}{0pt}{}{\lambda \in \sigma(T): |\Re (\lambda)| >1/t}{ \ or \ |\Im (\lambda)| >1/t}} \lambda\right|
&\le \sum_{\genfrac{}{}{0pt}{}{\lambda \in \sigma(T):}{ 1/t \le |\lambda| \le 2/t}} |\lambda| \\
&\le \frac2t \sum_{\genfrac{}{}{0pt}{}{\lambda \in \sigma(T):}{ 1/t \le |\lambda| \le 2/t}} 1 \\
&\le \frac2t d_T(\frac1t).
\end{align*}
Hence, by Lemma~\ref{l1}, we have 
$$\left|\omega\left(\frac1{\log(1+t)} \sum_{\lambda \in \sigma(T): |\lambda| >1/t} \lambda\right)-
\omega\left(\frac1{\log(1+t)} \sum_{\genfrac{}{}{0pt}{}{\lambda \in \sigma(T): |\Re (\lambda)| >1/t}{ \ or \ |\Im (\lambda)| >1/t}} \lambda\right)\right|$$
$$\le \omega\left(\frac2{t\log(1+t)} d_T(\frac1t)\right)=0.$$
\end{proof}

\begin{lem}\label{l3}
 For every normal operator $T\in \mathcal M_{1,\infty}$ and for every exponentiation invariant generalized limit $\omega$, we have
$$\omega\left(\frac1{\log(1+t)} \sum_{\genfrac{}{}{0pt}{}{\lambda \in \sigma(T): |\Re (\lambda)| >1/t}{ \ or \ |\Im (\lambda)| >1/t}} \Re (\lambda)\right)=
\omega\left(\frac1{\log(1+t)} \sum_{\lambda \in \sigma(T): |\Re (\lambda)| >1/t} \Re(\lambda)\right)
$$
and
$$\omega\left(\frac1{\log(1+t)} \sum_{\genfrac{}{}{0pt}{}{\lambda \in \sigma(T): |\Re (\lambda)| >1/t}{ \ or \ |\Im (\lambda)| >1/t}} \Im (\lambda)\right)=
\omega\left(\frac1{\log(1+t)} \sum_{\lambda \in \sigma(T): |\Im (\lambda)| >1/t} \Im (\lambda)\right). 
$$
\end{lem}
\begin{proof}
Consider the difference
\begin{align*}
\left|\sum_{\lambda \in \sigma(T): |\Re (\lambda)| >1/t} \Re(\lambda) - 
\sum_{\genfrac{}{}{0pt}{}{\lambda \in \sigma(T): |\Re (\lambda)| >1/t}{\ or \ |\Im (\lambda)| >1/t}} \Re(\lambda) \right|&
\le \sum_{\genfrac{}{}{0pt}{}{\lambda \in \sigma(T): |\Re (\lambda)| \le 1/t}{ \ and \ |\Im (\lambda)| >1/t}} |\Re(\lambda)| \\
&\le \frac1t \sum_{\lambda \in \sigma(T): |\Im (\lambda)| >1/t} 1.
\end{align*}

Since the operator $T$ is normal, it follows that $\Im (\sigma(T))= \sigma(\Im (T))$. Then 
$$\sum_{\lambda \in \sigma(T): |\Im (\lambda)| >1/t} 1 = d_{\Im (T)} (\frac1t)$$ 
and
$$\frac1{\log(1+t)}\left|\sum_{\lambda \in \sigma(T): |\Re (\lambda)| >1/t} \Re(\lambda) - 
\sum_{\genfrac{}{}{0pt}{}{\lambda \in \sigma(T): |\Re (\lambda)| >1/t}{ \ or \ |\Im (\lambda)| >1/t}} \Re(\lambda) \right|
\le \frac1{t \log(1+t)}d_{\Im (T)} (\frac1t).$$

Applying an exponentiation invariant generalized limit $\omega$ to both sides of the latter expression and using Lemma~\ref{l1}, we obtain the required equality.

The proof of the second equality is similar and is therefore omitted.
\end{proof}

The following theorem is an analogue of Lidskii formula for Dixmier traces generated by exponentiation invariant generalised limits.
\begin{thm}\label{LFP}
For every operator $T\in \mathcal M_{1,\infty}$ and for every exponentiation invariant generalized limit $\omega$, we have
$${\rm Tr}_\omega(T)=\omega\left(\frac1{\log(1+t)} \sum_{\lambda \in \sigma(T): |\lambda| >1/t} \lambda\right).$$
\end{thm}

\begin{proof}
For every compact operator $T$ there exist a compact normal operator $S$ 
and a compact quasi-nilpotent operator $Q$ such that $T=S+Q$ and $\sigma(T)=\sigma(S)$ (multiplicities coincide as well)~\cite[Theorems 1,6,7]{R}.
By the Weil theorem (see e.g.~\cite[Theorem 3.1]{GK}), the sequence $|\lambda(T)|$ is majorized by the sequence $\mu(T)$.
So, for $T\in \mathcal M_{1,\infty}$, we have $S,Q\in \mathcal M_{1,\infty}$.
By~\cite{K}, we have ${\rm Tr}_\omega(Q)=0$ for every quasi-nilpotent operator $Q$ and for every Dixmier trace.
Hence, it is sufficient to prove the statement of the theorem for a normal operator $T\in \mathcal M_{1,\infty}$.

Let $T\in \mathcal M_{1,\infty}$ be normal and let $\omega$ be an exponentiation invariant generalized limit.
Since $T$ is normal, it follows that $\Re (\sigma(T))= \sigma(\Re (T))$ and
$$\sum_{\lambda \in \sigma(T): |\Re (\lambda)| > 1/t } \Re(\lambda) = \sum_{\lambda \in \sigma(\Re(T)): |\lambda| > 1/t } \lambda.$$

Then by Corollary~\ref{c1} we obtain
$${\rm Tr}_\omega(\Re(T))= \omega\left(\frac1{\log(1+t)}\sum_{\lambda \in \sigma(T): |\Re (\lambda)| > 1/t } \Re(\lambda)\right).$$

Similarly, for the operator $\Im(T)$ we obtain
$${\rm Tr}_\omega(\Im(T))= \omega\left(\frac1{\log(1+t)}\sum_{\lambda \in \sigma(T): |\Im (\lambda)| > 1/t } \Im(\lambda)\right).$$ 

The assertion of the theorem follows from Lemma~\ref{l3} and Lemma~\ref{l2}.
\end{proof}

We complete this paper with an application to $\zeta$-functions of non-commutative geometry. 

Let $\omega$ be a generalized limit on $L_{\infty}$.
The functional $\zeta_{\omega}$ defined by the formula
$$\zeta_{\omega}(T)=\omega\left(\frac1r{\rm Tr}(T^{1+1/r})\right), \quad 0 \le T \in \mathcal{M}_{1,\infty}$$
is called a $\zeta$-function residue.
It is additive on the positive cone of $\mathcal{M}_{1,\infty}$ and, therefore, it extends to a fully symmetric functional on $\mathcal{M}_{1,\infty}$~\cite[Theorem 8]{SZ}.

The following theorem generalizes~\cite[Theorem 15]{SZ} by weakening the assumptions on a generalized limit $\omega$.

\begin{thm}\label{zeta} If $\omega$ is an exponentially invariant generalized limit, then 
$${\rm Tr}_\omega = \zeta_{\omega \circ \log}.$$
\end{thm}
\begin{proof}
It is sufficient to prove the equality ${\rm Tr}_\omega = \zeta_{\omega \circ \log}$ on positive operators $T \in \mathcal M_{1,\infty}$.

Define the function $\beta \ : \ (0,\infty) \to (0,\infty)$ by 
$$\beta(t) := \int_{e^{-t}}^\infty s \ d d_T(s)$$

and let
$$h(r):= \int_0^\infty e^{-t/r} \ d \beta(t).$$

It follows from spectral theorem that
$${\rm Tr} (T^{1+1/r})=h(r)$$

and by the definition of $\zeta$-function we have
$$\zeta_{\omega \circ \log}(T) = (\omega \circ \log)\left( \frac1r {\rm Tr} (T^{1+1/r})\right) = (\omega \circ \log)\left( \frac{h(r)}r\right).$$

Since $\omega \circ \log$ is a dilation invariant generalized limit, it follows from weak*-Karamata theorem (see e.g.~\cite[Theorem 2.2]{CPS2}) that
$$(\omega \circ \log)\left( \frac{h(r)}r\right) =(\omega \circ \log)\left( \frac{\beta(t)}t\right).$$

Hence, we obtain
\begin{align*}
\zeta_{\omega \circ \log}(T) &= (\omega \circ \log)\left( \frac{\beta(t)}t\right)\\ 
&= (\omega \circ \log)\left( \frac1t \int_{e^{-t}}^\infty s \ d d_T(s)\right) \\
&= \omega \left( \frac1{\log t} \int_{1/t}^\infty s \ d d_T(s)\right)\\
&={\rm Tr}_\omega(T),
\end{align*}
where the last equality is due to Theorem~\ref{LFP}.
\end{proof}

Research supported by the Australian Research Council.

\renewcommand{\refname}{References}


\begin{thebibliography}{X}
 \bibitem{BF} M. Benameur and T. Fack, Type II noncommutative geometry. I. Dixmier trace in von Neumann algebras, {\it Adv.Math.} {\bf 199} (2006), 29-87.
%
 \bibitem{CPRS1} A. Carey, J. Phillips, A. Rennie and F. Sukochev, The Hochschild class of the Chern character for semifinite spectral triples,  {\it J. Funct. Anal.} (1) {\bf213} (2004), 111--153.
%
 \bibitem{CPS2} 
  A. Carey, J. Phillips and F. Sukochev, Spectral flow and Dixmier traces, {\it Adv.Math.} (1) {\bf 173} (2003), 68--113.
%
  \bibitem{CS}  A. Carey and F. Sukochev, Dixmier traces and some application to the noncommutative geometry, (Russian) \textit{Uspehi. Mat. Nauk} (6) {\bf61} (2006), 45-110.  English translation in  \textit{Russian Math. Surveys}  (6) {\bf61} (2006), 1039--1099.
%
 \bibitem{CRSS} A. Carey, A. Rennie, A. Sedaev and F. Sukochev, The Dixmier trace and asymptotics of zeta functions, {\it J.Funct.Anal.} (2) {\bf 249} (2007), 253--283.
%
 \bibitem{C} A. Connes, {\it Noncommutative Geometry}, Academic Press, San Diego, 1994.
%
 \bibitem{D} J. Dixmier, Existence de traces non normales, {\it C. R. Acad. Sci. Paris} {\bf 262} (1966), A1107-A1108.

\bibitem{FK} T. Fack and H. Kosaki, Generalised $s$-numbers of
$\tau$-measurable operators, \textit{Pacific J. Math.} {\bf 123} (1986), 269--300.

\bibitem{GK} I. C. Gohberg and M. G. Krein, {\it Introduction to the Theory
of Non-selfadjoint Operators}, Translations of Mathematical Monographs,
{\bf 18}, AMS, 1969.

%
  \bibitem{KSS}
 N.J. Kalton, A. Sedaev and F. Sukochev, Fully symmetric functionals on a Marcinkiewicz space are Dixmier Traces, \textit{Adv. Math.} (4) {\bf 226} (2011), 3540--3549.
%
\bibitem{K} 
Kalton N. Spectral characterization of sums of commutators. I, \textit{J. Reine Angew. Math.} \textbf{504} (1998), 115--125.

 \bibitem{KPS}
 S. Krein, Ju. Petunin and E. Semenov, {\it Interpolation of linear operators,} Nauka, Moscow, 1978 (Russian). English translation in Translation of Mathematical Monographs, Amer. Math. Soc. {\bf 54} (1982).

%
  \bibitem{LSS}
 S. Lord, A. Sedaev, F. Sukochev, Dixmier traces as singular symmetric functionals and applications to measurable operators, {\it J. Funct. Anal.} (1) {\bf 224} (2005), 72--106.

 \bibitem{LP}
 W. A. J. Luxemburg and B. de Pagter,  Invariant Means for Positive Operators and Semigroups,
 {\it Katholieke Universiteit Nijmegen}, Subfaculteit Wiskunde, Nijmegen, (2005) 31--55.

\bibitem{R} 
J. Ringrose,  Super-diagonal forms for compact linear operators,  {\it Proc. London Math. Soc.} (3) \textbf{12} (1962) 367--384.
%
\bibitem{Sed}
A. Sedaev, Geometrical and topological aspects of interpolation spaces of Petre's K-method, Thesis, Voronezh 2010 (in Russian).

 \bibitem{SSZ}
 A. Sedaev, F. Sukochev and D. Zanin, Lidskii-type for\-mulae for Dixmier tra\-ces, {\it Int.Eq.Oper.Th.} (4) \textbf{68} (2010), 551--572.

 \bibitem{SUZ1}
 F. Sukochev, A. Usachev and D. Zanin, Generalized limits with additional invariance properties and their applications to noncommutative geometry, submitted manuscript.

 \bibitem{SUZ2}
 F. Sukochev, A. Usachev and D. Zanin, On the distinction between the classes of Dixmier and Connes-Dixmier traces, submitted manuscript.

\bibitem{SZ}
F. Sukochev, D. Zanin, {\it $\zeta-$function and heat kernel formulae}, J. Funct. Anal. {\bf 260} (2011), no. 8, 2451--2482.
\end{thebibliography}
\end{document}